\documentclass[12pt]{amsart}
\usepackage{amssymb}

\newtheorem{theorem}{Theorem}[section]
\newtheorem{definition}[theorem]{Definition}
\newtheorem{conjecture}[theorem]{Conjecture}
\newtheorem{proposition}[theorem]{Proposition}
\newtheorem{corollary}[theorem]{Corollary}

\newtheorem{lemma}[theorem]{Lemma}

\begin{document}

\title{Fusion rules for quantum reflection groups}

\author{Teodor Banica}
\address{T.B.: Department of Mathematics, Toulouse 3 University, 118 route de Narbonne, 31062 Toulouse, France. {\tt banica@math.ups-tlse.fr}}

\author{Roland Vergnioux} 
\address{R.V.: Department of Mathematics, Caen University, BP 5186, 14032 Caen Cedex, France. {\tt vergnioux@math.unicaen.fr}}


\begin{abstract}
We find the fusion rules for the quantum analogues of the complex reflection groups $H_n^s=\mathbb Z_s\wr S_n$. The irreducible representations can be indexed by the elements of the free monoid $\mathbb N^{*s}$, and their tensor products are given by formulae which remind the Clebsch-Gordan rules (which appear at $s=1$).
\end{abstract}

\maketitle

\section*{Introduction}

The last two decades have seen a remarkable unification process in mathematics and physics, with quantum groups playing a prominent role. The original discovery of Drinfeld \cite{dri} and Jimbo \cite{jim} was that the enveloping algebra $U(g)$ of a classical Lie algebra has a non-trivial deformation $U_q(g)$ in the category of Hopf algebras. Here the parameter $q$ can be complex, or formal.

Short after the constructions of Drinfeld and Jimbo, Woronowicz developed a general theory of compact quantum groups \cite{wo1}, \cite{wo2}. The algebras $U_q(g)$ with $q>0$ correspond to compact quantum groups, as shown by Rosso in \cite{ros}. 

Most of the study of compact quantum groups has focused on the extension and application of various differential geometry techniques. The work here, heavily influenced by Connes' book \cite{con}, currently follows a number of independent directions, belonging to noncommutative geometry. See \cite{kye}, \cite{ntu}.

Another part of work has gone into the study of free quantum groups. These were introduced in two papers of Wang \cite{wa1}, \cite{wa2}. The idea is as follows: let $G\subset U_n$ be a compact group. The $n^2$ matrix coordinates $u_{ij}$ satisfy certain relations $R$, and generate the algebra $C(G)$. One can define then the universal algebra $A=C^*(u_{ij}|R)$, and for a suitable choice of the relations $R$ we get a Hopf algebra. We have the heuristic formula $A=C(G^+)$, where $G^+$ is a compact quantum group, called free version of $G$. Observe that we have $G\subset G^+$.

This construction is not axiomatized, in the sense that $G^+$ depends on the relations $R$, and it is not known in general what the good choice of $R$ is. For instance any choice with $R$ including the commutativity relations $u_{ij}u_{kl}=u_{kl}u_{ij}$ would be definitely a bad one, because in this case we would get $G^+=G$. Moreover, any choice with $R$ including certain relations which imply these commutativity relations would be a bad one as well, and the problem comes from here.

There are a number of general methods, however, which can be used in order to understand the correct formulation of the operation $G\to G^+$. One of them, heavily used in the context of quantum permutation groups \cite{bbi}, \cite{bbc} consists in saying that ``$G^+$ should be the quantum symmetry group of the object $X$ whose symmetry group is $G$''. The axiomatization work here seems to evolve towards the quite general situation where $X$ is a spectral triple in the sense of Connes \cite{con}, thanks to the recent work of Goswami \cite{gos}.

For the known examples of free quantum groups, the main problem is to compute certain representation theory invariants. As in the case of classical groups, the central problem is that of classifying the irreducible representations, and finding their fusion rules. The story here is as follows:
\begin{enumerate}
\item The first two examples of free quantum groups are $O_n^+$ and $U_n^+$. These quantum groups were constructed in Wang's thesis \cite{wa1}, and their fusion rules were found in the first author's thesis \cite{ba1}.

\item The knowledge of the fusion rules for $O_n^+$ and $U_n^+$ allows the study of their duals, and some preliminary results were obtained in \cite{ba1}. The systematic study in this sense has begun with the second author's thesis \cite{ve1}.

\item The third free quantum group is $S_n^+$, constructed in Wang's paper \cite{wa2}. The fusion rules for $S_n^+$ were found in \cite{ba2}, with the quite surprising result that these are the same as those for $SO_3$, independently of $n\geq 4$.

\item The analytic study of the duals of $O_n^+$ and $U_n^+$ has been intensively developed in the last 5 years, with a number of key results obtained by the second author, Vaes, and their collaborators \cite{vvv}, \cite{vve}, \cite{ve2}.

\item In the meantime, the construction of new free quantum groups, along with some preliminary axiomatization work, has been pursued by the first author, Bichon, Collins, and their collaborators \cite{bb+}, \cite{bbi}, \cite{bbc}. 
\end{enumerate}

The purpose of the present work is to make a bridge between (4) and (5), with an explicit computation of fusion rules, in the spirit of \cite{ba1}, \cite{ba2}. It is our hope that, as it was the case with \cite{ba1}, \cite{ba2}, the present results will be of help for both approaches. Some comments in this sense will be given in the end of the paper.

The main object of interest will be the complex reflection group $H_n^s$. This is the group of monomial $n\times n$ matrices (i.e. one nonzero entry on each row and each column), having as nonzero entries the $s$-th roots of unity. 

The construction of $H_n^{s+}$ was done in several steps, by overcoming a series of quite unexpected algebraic obstacles, the story being as follows:
\begin{enumerate}
\item At $s=1$ we have the symmetric group $S_n$, investigated by Wang in \cite{wa2}. The subtlety here comes from the fact that $S_n^+$ is not a finite quantum group. However, $S_n^+$ makes sense as a compact quantum group.

\item At $s=2$ we have the hyperoctahedral group $H_n$, studied in \cite{bbc}. The surprise here comes from the fact that $H_n^+$ is not the quantum symmetry group of the hypercube (which can be shown to be $O_n^{-1}$).

\item The general case $s\geq 1$ was worked out in \cite{bb+}. Once again, there were several candidates here for $H_n^{s+}$, and a quite lenghty free probability computation had to be done in order to find the good one.

\item The general series of complex reflection groups is $\{H_n^{rs}\}$, where $H_n^{rs}\subset H_n^s$ with $r|s$ consists of matrices whose product of nonzero entries is a $r$-th root of unity. The existence of $H_n^{rs+}$ is an open, uncertain problem.
\end{enumerate}

In this paper we classify the irreducible representations of $H_n^{s+}$, and we find their fusion rules. The main result states that the irreducible representations can be indexed by the words over $\mathbb Z_s$, and the fusion rules are of the following type:
$$r_x\otimes r_y=\sum_{x=vz,y=\bar{z}w}r_{vw}+r_{v\cdot w}$$

We refer to the section 7 below for the precise meaning of this formula, and to section 8 below for some equivalent formulations.

The above formula can be thought of as being a ``level $s$'' generalization of the Clebsch-Gordan formula for $SO_3$, which appears at $s=1$.

The proof is considerably more complicated than the one in \cite{ba2} at $s=1$, and uses various techniques from \cite{bbi}, \cite{bbc}, and from the main paper \cite{bb+}.

The paper is organized as follows. In 1-3 we discuss the construction and the basic algebraic properties of the quantum reflection groups. In 4-7 we work out the fusion rules, and in 8-9 we discuss some related questions.

The final section, 10, contains a conjectural statement regarding the fusion rules of arbitrary free quantum groups, along with a few concluding remarks.

\section{Quantum reflection groups}

A square matrix is called monomial if it has exactly one nonzero entry on each row and on each column. The basic examples are the permutation matrices.

\begin{definition}
$H_n^s$ is the group of monomial $n\times n$ matrices having as nonzero entries the $s$-th roots of unity.
\end{definition}

In other words, an element of $H_n^s$ is a permutation matrix, with the 1 entries replaced by certain $s$-th roots of unity. Observe that we have $H_n^s=\mathbb Z_s\wr S_n$.

We allow the value $s=\infty$ in the above definition, with the convention that a root of unity of infinite order is nothing but an element on the unit circle.

Of special interest among the groups $H_n^s$ are those corresponding to the values $s=1,2,\infty$. These groups will be used as key examples for all the considerations in this paper, and most definitions or theorems will be illustrated in this way.

\begin{proposition}
The groups $H_n^s$ with $s=1,2,\infty$ are as follows.
\begin{enumerate}
\item $H_n^1=S_n$ is the symmetric group.

\item $H_n^2=H_n$ is the hyperoctahedral group.

\item $H_n^\infty=K_n$ is the group of unitary monomial matrices.
\end{enumerate}
\end{proposition}

The groups $H_n^s$ are in fact part of a more general series. Indeed, for $r|s$ we can consider the subgroup $H_n^{rs}\subset H_n^s$ formed by the matrices whose product of nonzero entries is a $r$-th root of unity. Observe that we have $H_n^{ss}=H_n^s$. 

The groups $H_n^{rs}$ form the general series of complex reflection groups. We should probably mention that the standard group theory notation for this series is $G(d,de,n)=H_n^{d,de}$, so in particular we have $H_n^s=G(s,s,n)$.

The following key definition is from \cite{bb+}.

\begin{definition}
$A_h^s(n)$ is the universal $C^*$-algebra generated by $n^2$ normal elements $u_{ij}$, subject to the following relations $R$:
\begin{enumerate}
\item $u=(u_{ij})$ is unitary.

\item $u^t=(u_{ji})$ is unitary.

\item $p_{ij}=u_{ij}u_{ij}^*$ is a projection.

\item $u_{ij}^s=p_{ij}$.
\end{enumerate}
\end{definition}

We use here the standard operator algebra terminology: an element $a$ is called normal if $aa^*=a^*a$, unitary if $aa^*=a^*a=1$, and projection if $a^2=a^*=a$.

We allow the value $s=\infty$ in the above definition, with the convention that the last axiom simply disappears in this case.

Observe that for $s<\infty$ the normality condition is actually redundant. This is because a partial isometry $a$ subject to the relation $aa^*=a^s$ is normal.

It follows from definitions, and from standard operator algebra tricks, that $A_h^s(n)$ is a Hopf algebra, with comultiplication, counit and antipode as follows:
\begin{eqnarray*}
\Delta(u_{ij})&=&\sum_{k=1}^n u_{ik}\otimes u_{kj}\\
\varepsilon(u_{ij})&=&\delta_{ij}\\
S(u_{ij})&=&u_{ji}^*
\end{eqnarray*}

More precisely, the above formulae show that $A_h^s(n)$ is a Hopf $C^*$-algebra, in the sense of Woronowicz's fundamental paper \cite{wo1}.

The relation with the group $H_n^s$ is as follows.

\begin{proposition}
$C(H_n^s)$ is isomorphic to the universal commutative $C^*$-algebra generated by $n^2$ abstract variables $u_{ij}$, subject to the relations $R$.
\end{proposition}

This assertion follows indeed from the Stone-Weierstrass theorem and from the Gelfand theorem. See \cite{bb+}.

The comparison between Definition 1.3 and Proposition 1.4 leads to the conclusion that we have a surjective morphism of $C^*$-algebras $\pi:A_h^s(n)\to C(H_n^s)$, whose kernel is the commutator ideal of $A_h^s(n)$.

The best interpretation, however, is in terms of quantum groups. It follows from definitions that $\pi$ is a Hopf algebra morphism, so if $H_n^{s+}$ is the compact quantum group associated to $A_h^s(n)$, then we have an embedding $H_n^s\subset H_n^{s+}$. 

\begin{proposition}
The algebras $A_h^s(n)$ with $s=1,2,\infty$, and their presentation relations in terms of the entries of the matrix $u=(u_{ij})$, are as follows.
\begin{enumerate}
\item For $A_h^1(n)=A_s(n)$, the matrix $u$ is magic: all its entries are projections, summing up to $1$ on each row and column.

\item For $A_h^2(n)=A_h(n)$ the matrix $u$ is cubic: it is orthogonal, and the products of pairs of distinct entries on the same row or the same column vanish.

\item For $A_h^\infty(n)=A_k(n)$ the matrix $u$ is unitary, its transpose is unitary, and all its entries are normal partial isometries.
\end{enumerate}
\end{proposition}

In this statement (1) and (2) follow from definitions and from standard operator algebra tricks, and (3) is just a translation of the definition (a partial isometry is an element $a$ such that $aa^*$ and $a^*a$ are projections). See \cite{bb+}.

\section{Finiteness considerations}

In this section and in the next one we find the quantum analogues of some basic properties of $H_n^s$. We should mention that these results were proved in \cite{bbc} in the case $s=2$, and were announced in \cite{bb+} in the general case $s\geq 1$.

The most obvious property of $H_n^s$ is that this is a finite group. 

The quantum analogue of this result, however, doesn't hold, even at $s=1$, as shown by the following result of Wang \cite{wa2}.

\begin{theorem}
For $n=2,3$ the canonical map $A_s(n)\to C(S_n)$ is an isomorphism. For $n\geq 4$ the algebra  $A_s(n)$ is not commutative, and infinite dimensional. 
\end{theorem}

In other words, the compact quantum group $S_n^+$ with $n\geq 4$ is a not a classical group, nor a finite quantum group.

Now back to the general case $s\geq 1$, the correct approach to this quite mysterious finiteness issue is as follows:
\begin{enumerate}
\item First, we will prove that $H_n^{s+}$ is a quantum permutation group, in analogy with the fact that $H_n^s$ is a permutation group. 

\item Then we will prove that we have a free wreath product decomposition $H_n^{s+}=\mathbb Z_s\wr_*S_n^+$, in analogy with the decomposition $H_n^s=\mathbb Z_s\wr S_n$.
\end{enumerate}

In Hopf algebra terms, we first have to prove that $A_h^s(n)$ is a quantum permutation algebra, i.e. is a quotient of a suitable Wang algebra $A_s(N)$. For this purpose, we use the method of ``sudoku matrices'' from \cite{bbc}.

We recall from Proposition 1.5 that a magic unitary is a square matrix of projections, which sum up to $1$ on each row and column.

\begin{definition}
A $(s,n)$-sudoku matrix is a magic unitary of size $sn$, of the form
$$m=\begin{pmatrix}
a^0&a^1&\ldots&a^{s-1}\\
a^{s-1}&a^0&\ldots&a^{s-2}\\
\ldots&\ldots&\ldots&\ldots\\
a^1&a^2&\ldots&a^0
\end{pmatrix}$$
where $a^0,\ldots,a^{s-1}$ are $n\times n$ matrices.
\end{definition}

Observe that $m$ is a circulant matrix. By using some standard tensor product identifications, and modulo $s$ indices, we can use the following writing for it:
$$m=\left(a^{q-p}_{ij}\right)_{pi,qj}$$

The basic example of such matrices is in relation with the group $H_n^s$. With the notation $w=e^{2\pi i/s}$, each of the $n^2$ matrix coordinates $u_{ij}:H_n^s\to\mathbb C$ takes values in the set $\{0\}\cup\{1,w,\ldots,w^{s-1}\}$, hence decomposes as follows:
$$u_{ij}=\sum_{r=0}^{s-1}w^ra^r_{ij}$$

Here each $a^r_{ij}$ is by definition a function taking values in $\{0,1\}$. We see that each $a^r_{ij}$ is a projection in the $C^*$-algebra sense, and it follows from definitions that these projections form a sudoku matrix, in the above sense.

\begin{theorem}
We have the following results.
\begin{enumerate}
\item $C(H_n^s)$ is isomorphic to the universal commutative $C^*$-algebra generated by the entries of a $(s,n)$-sudoku matrix.

\item $A_h^s(n)$ is isomorphic to the universal $C^*$-algebra generated by the entries of a $(s,n)$-sudoku matrix.
\end{enumerate}
\end{theorem}

\begin{proof}
(1) This assertion, which is included here for symmetry reasons, and which won't be used in what follows, can be proved either directly, or by using (2) and the fact that $C(H_n^s)$ is the maximal commutative quotient of $A_h^s(n)$. 

(2) We denote by $A$ the universal algebra in the statement. According to Definition 2.1, we have the following presentation formula:
$$A=C^*\left(a_{ij}^p\ \Big\vert \left(a^{q-p}_{ij}\right)_{pi,qj}=(s,n)-\mbox{ sudoku }\right)$$

Consider also the algebra $A_h^s(n)$. According to Definition 1.3, this is presented by certain relations $R$, that we call here level $s$ cubic conditions:
$$A_h^s(n)=C^*\left(u_{ij}\ \Big\vert\  u=n\times n\mbox{ level $s$ cubic }\right)$$

We will construct a pair of inverse morphisms between these algebras.

{\bf Step 1.} Consider the following matrix:
$$U_{ij}=\sum_pw^{-p}a^p_{ij}$$ 

We claim that this is a level $s$ cubic unitary. Indeed, by using the sudoku condition, the verification of (1-4) in Definition 1.3 goes as follows.

(1) The fact that we have $UU^*=1$ can be checked as follows:
\begin{eqnarray*}
(UU^*)_{ij}
&=&\sum_{kpq}w^{-p}a^p_{ik}w^qa^q_{jk}\\
&=&\sum_{pq}w^{q-p}\sum_ka^p_{ik}a^q_{jk}\\
&=&\sum_{pq}w^{q-p}\delta_{pq}\delta_{ij}\sum_ka^p_{ik}\\
&=&\delta_{ij}\sum_{pk}a^p_{ik}\\
&=&\delta_{ij}
\end{eqnarray*}

By symmetry reasons, the verification of $U^*U=1$ is similar.

(2) The verification of $U^t\bar{U}=1$ and $\bar{U}U^t=1$ is similar.

(3) We first compute the elements $P_{ij}=U_{ij}U_{ij}^*$:
\begin{eqnarray*}
P_{ij}
&=&\sum_{pq}w^{-p}a^p_{ij}w^qa^q_{ij}\\
&=&\sum_{pq}w^{q-p}a^p_{ij}a^q_{ij}\\
&=&\sum_pa^p_{ij}
\end{eqnarray*}

This is a sum of pairwise orthogonal projections, so it is a projection.

(4) We compute now the $s$-th power of $U_{ij}$:
\begin{eqnarray*}
U_{ij}^s
&=&\left(\sum_pw^{-p}a_{ij}^p\right)^s\\
&=&\sum_p\left(w^{-p}a^p_{ij}\right)^s\\
&=&\sum_pa^p_{ij}\\
&=&P_{ij}
\end{eqnarray*}

Summarizing, the elements $U_{ij}$ form a level $s$ cubic matrix, so we can define a morphism $\Phi:A_h^s(n)\to A$ by the formula $\Phi(u_{ij})=U_{ij}$.

{\bf Step 2.} Consider the following elements, with the convention $u_{ij}^0=p_{ij}$:
$$A^p_{ij}=\frac{1}{s}\sum_rw^{rp}u^r_{ij}$$

It follows from the cubic condition that these elements form a level $s$ sudoku unitary, with the verification going as follows:

(1) First, these elements are self-adjoint:
\begin{eqnarray*}
(A^p_{ij})^*
&=&\frac{1}{s}\sum_rw^{-rp}(u^r_{ij})^*\\
&=&\frac{1}{s}\sum_rw^{-rp}u^{s-r}_{ij}\\
&=&\frac{1}{s}\sum_rw^{(s-r)p}u^{s-r}_{ij}\\
&=&A^p_{ij}
\end{eqnarray*}

(2) We check now that these elements are idempotents:
\begin{eqnarray*}
(A^p_{ij})^2
&=&\frac{1}{s^2}\sum_{rt}w^{rp}u^r_{ij}w^{tp}u^t_{ij}\\
&=&\frac{1}{s^2}\sum_{rt}w^{(r+t)p}u_{ij}^{r+t}\\
&=&\frac{1}{s}\sum_lw^{lp}u_{ij}^l\\
&=&A_{ij}^p
\end{eqnarray*}

(3) We compute the sum on the rows of $M=(A^{q-p}_{ij})_{pi,qj}$:
\begin{eqnarray*}
\sum_{jp}A^p_{ij}
&=&\frac{1}{s}\sum_{jpr}w^{rp}u_{ij}^r\\
&=&\frac{1}{s}\sum_{jr}u_{ij}^r\sum_pw^{rp}\\
&=&\sum_ju_{ij}^0\\
&=&1
\end{eqnarray*}

(4) By symmetry reasons, the sum on the columns of $M$ is $1$ as well.

Summarizing, the elements $A^p_{ij}$ form a sudoku unitary, so we can define a morphism $\Psi:A\to A_h^s(n)$ by the formula $\Psi(a^p_{ij})=A^p_{ij}$.

{\bf Step 3.} We check now the fact that $\Phi,\Psi$ are indeed inverse morphisms:
\begin{eqnarray*}
\Psi\Phi(u_{ij})
&=&\sum_pw^{-p}A^p_{ij}\\
&=&\frac{1}{s}\sum_pw^{-p}\sum_rw^{rp}u_{ij}^r\\
&=&\frac{1}{s}\sum_{pr}w^{(r-1)p}u_{ij}^r\\
&=&u_{ij}
\end{eqnarray*}

As for the other composition, we have the following computation:
\begin{eqnarray*}
\Phi\Psi(a^p_{ij})
&=&\frac{1}{s}\sum_rw^{rp}U_{ij}^r\\
&=&\frac{1}{s}\sum_rw^{rp}\sum_qw^{-rq}a_{ij}^q\\
&=&\frac{1}{s}\sum_qa_{ij}^q\sum_rw^{r(p-q)}\\
&=&a^p_{ij}
\end{eqnarray*}

This finishes the proof.
\end{proof}

\section{Algebraic structure}

We know from the previous section that $H_n^{s+}$ is a quantum permutation group, in analogy with the fact that $H_n^s$ is a permutation group. In this section we discuss the quantum analogue of the decomposition $H_n^s=\mathbb Z_s\wr S_n$. 

\begin{lemma}
A $sn\times sn$ magic unitary commutes with the square matrix
$$\Sigma=
\begin{pmatrix}
0&I_n&0&\ldots&0\cr
0&0&I_n&\ldots&0\cr
\ldots&\ldots&\ldots&\ldots&\ldots\cr
0&0&0&\ldots&I_n\cr
I_n&0&0&\ldots&0
\end{pmatrix}$$
if and only if it is a sudoku matrix in the sense of Definition 2.2.
\end{lemma}

\begin{proof}
The commutation with $\Sigma$ ensures indeed that the matrix is circulant.
\end{proof}

Let $C_s$ be the oriented cycle with $s$ vertices, and consider the graph $C_s^n$ consisting of $n$ disjoint copies of it. Observe that, with a suitable labeling of the vertices, the adjacency matrix of this graph is the above matrix $\Sigma$.

The quantum symmetry algebra of a finite graph is the quotient of the quantum permutation algebra on the set of vertices by the relations making the fundamental corepresentation commute with the adjacency matrix. See \cite{bbi}.

\begin{theorem}
We have the following results.
\begin{enumerate}
\item $H_n^s$ is the symmetry group of $C_s^n$.

\item $A_h^s(n)$ is the quantum symmetry algebra of $C_s^n$.
\end{enumerate}
\end{theorem}

\begin{proof}
(1) follows from definitions, and (2) follows from Theorem 2.3 and Lemma 3.1. Indeed, $A_h^s(n)$ is the quotient of $A_s(sn)$ by the relations making the fundamental corepresentation commute with the adjacency matrix of $C_s^n$.
\end{proof}

According to the work of Bichon \cite{bic}, the free analogue of the notion of wreath product is that of free wreath product at the level of Hopf algebras. 

\begin{definition}
The free wreath product of two quantum permutation algebras $(A,u)$ and $(B,v)$ is given by
$$A*_wB=(A^{*n}*B)/<[u_{pq}^{(i)},v_{ij}]=0>$$
where $n$ is the size of $v$, with magic unitary matrix $w_{pi,qj}=u_{pq}^{(i)}v_{ij}$.
\end{definition}

This definition is justified by formulae of the following type, where $G,A$ denote classical symmetry groups, respectively quantum
symmetry algebras:
\begin{eqnarray*}
G(X*Y)&=&G(X)\ \wr\ G(Y)\\
A(X*Y)&=&A(X)*_wA(Y)
\end{eqnarray*}

There are several such formulae, depending on the types of graphs and products considered. See \cite{bic}, \cite{bbi}. The formula we are interested in is:
\begin{eqnarray*}
G(X\ldots X)&=&G(X)\ \wr \ G(\circ\ldots\circ)\\
A(X\ldots X)&=&A(X)*_wA(\circ\ldots\circ)\end{eqnarray*}

Here $X$ is a finite graph, $\circ$ is a point, and the dots mean $n$-fold disjoint union. For the precise statement and proof of this result, we refer to \cite{bbi}.

We are now in position of stating the main result in this section.

\begin{theorem}
We have the following results.
\begin{enumerate}
\item $H_n^s=\mathbb Z_s\wr S_n$.

\item $A_h^s(n)=C(\mathbb Z_s)*_wA_s(n)$.
\end{enumerate}
\end{theorem}

\begin{proof}
This follows from Theorem 3.2 and from the above formulae, first established in \cite{bic} and later on refined in \cite{bbi}, by using the graph $X=C_s$.

Observe that (1) is in fact clear from definitions. We would like to present below a self-contained proof of (2), by constructing a pair of inverse morphisms. This will compress the above-mentioned combined arguments from \cite{bic}, \cite{bbi}.

{\bf Step 1.} First, we have to fix some notations for the algebra on the right. We view $\mathbb Z_s$ as the group formed by the powers of the basic cyclic matrix:
$$\sigma=
\begin{pmatrix}
0&1&0&\ldots&0\cr
0&0&1&\ldots&0\cr
\ldots&\ldots&\ldots&\ldots&\ldots\cr
0&0&0&\ldots&1\cr
1&0&0&\ldots&0
\end{pmatrix}$$
 
Thus we have $\mathbb Z_s\subset M_s(\mathbb C)$, and the magic unitary $u$ corresponding to the quantum permutation algebra $C(\mathbb Z_s)$ is the matrix of coordinates on $\mathbb Z_s$:
$$u_{pq}(\sigma^r)=\delta_{q-p,r}$$

Here, and in what follows, all the indices $p,q,r,\ldots$ are taken mod $s$. Observe that $u$ is a circulant matrix, and in particular we have:
$$u_{pq}=u_{0,q-p}$$

{\bf Step 2.} We construct a morphism $\Phi:A_h^s(n)\to C(\mathbb Z_s)*_wA_s(n)$.

Consider the standard generators $w_{pi,qj}=u_{pq}^{(i)}v_{ij}$ of the algebra on the right, as in Definition 3.3. We claim that the following elements form a sudoku unitary:
$$A_{ij}^p=u_{0p}^{(i)}v_{ij}$$

Indeed, the corresponding matrix $M$ as in Definition 2.2 is given by:
\begin{eqnarray*}
M_{pi,qj}
&=&A_{ij}^{q-p}\\
&=&u_{0,q-p}^{(i)}v_{ij}\\
&=&u_{pq}^{(i)}v_{ij}\\
&=&w_{pi,qj}
\end{eqnarray*}

Since this latter matrix is known to be magic, the elements $A_{ij}^p$ form indeed a sudoku unitary, and we get a morphism $\Phi$ as claimed.

{\bf Step 3.} We construct a morphism $\Psi:C(\mathbb Z_s)*_wA_s(n)\to A_h^s(n)$.

Consider the standard sudoku generators $a_{ij}^p$ of the algebra on the right, as in Definition 2.2. We define elements $U_{pq}^{(i)}$ and $V_{ij}$ as follows:
$$U_{pq}^{(i)}=\sum_ka_{ik}^{q-p}$$
$$V_{ij}=\sum_ra_{ij}^r$$

It is routine to check that each of the matrices $U^{(i)}$ produces a morphism $C(\mathbb Z_s)\to A_h^s(n)$, and that $V$ produces a morphism $A_s(n)\to A_h^s(n)$. Moreover, we have the following commutation relation:
\begin{eqnarray*}
[U_{pq}^{(i)},V_{ij}]
&=&\left[\sum_ka_{ik}^{q-p},\sum_ra_{ij}^r\right]\\
&=&\sum_{kr}[a_{ik}^{q-p},a_{ij}^r]\\
&=&0
\end{eqnarray*}

Summarizing, the elements $U_{pq}^{(i)}$ and $V_{ij}$ satisfy the defining relations for the free wreath product, so we get a morphism $\Psi$ as claimed.

{\bf Step 4.} We check now that fact that $\Phi,\Psi$ are indeed inverse morphisms. In one sense, we have the following computation:
\begin{eqnarray*}
\Psi\Phi(a_{ij}^p)
&=&U_{0p}^{(i)}V_{ij}\\
&=&\sum_{kr}a_{ik}^pa_{ij}^r\\
&=&a_{ij}^p
\end{eqnarray*}

As for the other composition, we have the following computation:
\begin{eqnarray*}
\Phi\Psi(w_{pi,qj})
&=&\sum_{kr}A_{ik}^{q-p}A_{ij}^r\\
&=&\sum_{kr}u_{0,q-p}^{(i)}v_{ik}u_{0r}^{(i)}v_{ij}\\
&=&\sum_{kr}u_{0,q-p}^{(i)}u_{0r}^{(i)}v_{ik}v_{ij}\\
&=&u_{0,q-p}^{(i)}v_{ij}\\
&=&w_{pi,qj}
\end{eqnarray*}

This finishes the proof.
\end{proof}

We can use the above result in order to deduce the structure of $A_h^s(n)$ in the cases $n=2,3$, which are to be avoided in what follows.

\begin{corollary}
The algebras $A_h^s(n)$ with $n=2,3$ are as follows:
\begin{enumerate}
\item $A_h^s(2)=C(\mathbb Z_s)*_wC(\mathbb Z_2)$.

\item $A_h^s(3)=C(\mathbb Z_s)*_wC(S_3)$.
\end{enumerate}
\end{corollary}

\begin{proof}
This follows from Theorem 2.1 and Theorem 3.4.
\end{proof}

\section{Basic corepresentations}

In this section and in the next few ones we discuss the classification of irreducible corepresentations of $A_h^s(n)$, and the computation of their fusion rules.

We recall that, according to Woronowicz's fundamental paper \cite{wo1}, an analogue of the Peter-Weyl theory is available for the compact quantum groups. 

In Hopf algebra terms, the objects of interest are the finite dimensional irreducible unitary corepresentations, in the following sense. 

\begin{definition}
A finite dimensional unitary corepresentation of a Hopf $C^*$-algebra $A$ is a unitary matrix $u\in M_n(A)$ satisfying the following conditions:
\begin{eqnarray*}
\Delta(u_{ij})&=&\sum_{k=1}^nu_{ik}\otimes u_{kj}\\
\varepsilon(u_{ij})&=&\delta_{ij}\\
S(u_{ij})&=&u_{ji}^*
\end{eqnarray*}
Such a corepresentation is called irreducible if the matrices $T\in M_n(\mathbb C)$ commuting with it, $Tu=uT$, reduce to the scalar multiples of the identity.
\end{definition}

The sum and tensor product of two corepresentations $u,v$ are by definition the matrices $u+v={\rm diag}(u,v)$ and $u\otimes v=(u_{ij}v_{ab})_{ia,jb}$. 

The basic examples of corepresentations are the fundamental one $u=(u_{ij})$, and its complex conjugate $\bar{u}=(u_{ij}^*)$. In this section we use $u$ and $\bar{u}$ in order to construct a whole family of ``basic'' corepresentations of $A_h^s(n)$.

For this purpose, we go back to the elements $u_{ij},p_{ij}$ in Definition 1.3. We recall that, as a consequence of Proposition 1.5, $p$ is a magic unitary.

\begin{lemma}
The elements $u_{ij}$ and $p_{ij}$ satisfy:
\begin{enumerate}
\item $p_{ij}u_{ij}=u_{ij}$.

\item $u_{ij}^*=u_{ij}^{s-1}$.

\item $u_{ij}u_{ik}=0$ for $j\neq k$.
\end{enumerate}
\end{lemma}

\begin{proof}
We use the fact that in a $C^*$-algebra, $aa^*=0$ implies $a=0$.

(1) With $a=(p_{ij}-1)u_{ij}$ we have:
\begin{eqnarray*}
aa^*
&=&(p_{ij}-1)u_{ij}u_{ij}^*(p_{ij}-1)\\
&=&(p_{ij}-1)p_{ij}(p_{ij}-1)\\
&=&0
\end{eqnarray*}

Thus we have $a=0$, which gives the result.

(2) With $a=u_{ij}^*-u_{ij}^{s-1}$ we have $aa^*=0$, which gives the result.

(3) With $a=u_{ij}u_{ik}$ we have $aa^*=0$, which gives the result.
\end{proof}

In what follows, we make the convention $u_{ij}^0=p_{ij}$.

\begin{theorem}
The algebra $A_h^s(n)$ has a unique family of $n$-dimensional corepresentations $\{u_k|k\in\mathbb Z\}$, satisfying the following conditions:
\begin{enumerate}
\item $u_k=(u_{ij}^k)$ for any $k\geq 0$.

\item $u_k=u_{k+s}$ for any $k\in\mathbb Z$.

\item $\bar{u}_k=u_{-k}$ for any $k\in\mathbb Z$.
\end{enumerate}
\end{theorem}

\begin{proof}
We first prove that the matrix $u_k=(u_{ij}^k)$ is a corepresentation, for any $k\geq 1$. By using the last assertion in the previous lemma, we get:
\begin{eqnarray*}
\Delta(u_{ij}^k)
&=&(\Delta(u_{ij}))^k\\
&=&\left(\sum_l u_{il}\otimes u_{lj}\right)^k\\
&=&\sum_{l_1\ldots l_k}u_{il_1}\ldots u_{il_k}\otimes u_{l_1j}\ldots u_{l_kj}\\
&=&\sum_lu_{il}^k\otimes u_{lj}^k
\end{eqnarray*}

As for the formulae $\varepsilon(u_{ij}^k)=\delta_{ij}$ and $S(u_{ij}^k)=u_{ji}^{*k}$, these follow from $\varepsilon(u_{ij})=\delta_{ij}$ and $S(u_{ij})=u_{ji}^*$, by using the multiplicative properties of $\varepsilon,S$.

We claim now that we have $u_k=u_{k+s}$, for any $k\geq 1$. Indeed, this follows from the last assertion in the previous lemma:
\begin{eqnarray*}
u_{ij}^{k+s}
&=&u_{ij}^ku_{ij}^s\\
&=&u_{ij}^kp_{ij}\\
&=&u_{ij}^k
\end{eqnarray*}

Summarizing, the conditions (1) and (2) in the statement define a unique family of $n$-dimensional corepresentations $\{u_k|k\in\mathbb Z\}$, and it remains to check that these corepresentations satisfy the condition (3). But this latter condition follows from the second assertion in the previous lemma, and we are done.
\end{proof}

\section{Noncrossing partitions}

In this section and in the next one we compute the intertwiners between the various tensor products between the basic corepresentations $u_i$.

The idea is to use the canonical arrow $A_h^s(n)\to A_s(n)$. This maps all the corepresentations $u_i$ into $U$, the fundamental corepresentation of $A_s(n)$, so by functoriality we get embeddings as follows:
$$Hom(u_{i_1}\otimes\ldots\otimes u_{i_k},u_{j_1}\otimes\ldots\otimes u_{j_l})\subset Hom(U^{\otimes k},U^{\otimes l})$$

So, our first task will be to present a detailed description of the spaces on the right. Then, a careful study will allow us to find the spaces on the left.

Recall from \cite{spe} the following definition of noncrossing partitions of an ordered set $S$. A partition $S = P_1 \sqcup P_2 \sqcup \cdots \sqcup P_k$ is called noncrossing if, for any distinct classes $P_i = (s_1 < s_2 < \cdots < s_l)$ and $P_j = (t_1 < t_2 < \cdots < t_m)$ of the partition we have:
$$ t_k < s_1 < t_{k+1} ~~\Longleftrightarrow~~ t_k < s_l < t_{k+1} $$
Such a partition can be pictorially represented by putting the elements of $S$ on the real line and joining together the elements of each $P_i$ by strings in the upper half plane, in such a way that the strings of the resulting picture do not intersect --- see the example after the Definition.

\begin{definition}
We denote by $NC(k,l)$ the set of noncrossing partitions of the set with repetitions $\{1,\ldots,k,1,\ldots,l\}$ ordered as $1 < \cdots < k < l < \cdots < 1$. These will be pictured as
$$p=\left\{\begin{matrix}
1\ldots k\cr
\mathcal P\cr
1\ldots l
\end{matrix}\right\}$$
where $\mathcal P$ is a noncrossing diagram joining the elements in the same class of the partition.
\end{definition}

Observe that $NC(k,l)$ is in correspondence with the set $NC(k+l)$ of noncrossing partitions of $\{1,\ldots,k+l\}$. As an example, consider the following partition in $NC(6)$:
$$p=\{1,2,5\}\cup\{3,4\}\cup\{6\}$$

The corresponding element of $NC(6,0)$ is pictured as follows:
$$p_{60}=\left\{\begin{matrix}
1\ \ 2\  3\, 4\ 5\ \ 6\cr
|_{{\ }_{\!\!\!\!-\!\!-\!\!-\!\!-\!\!}}^
{\ }|_{{\ }_{\!\!\!\!-\!\!-\!\!-\!\!-\!\!-\!\!-\!\!-\!\!}}^{\ \ \sqcup}|\ \ \ |
\end{matrix}\right\}$$

The corresponding element of $NC(0,6)$ is pictured as follows:
$$p_{06}=\left\{\begin{matrix}|^{{\ }^{\!\!\!\!-\!\!-\!\!-\!\!-\!\!}}_
{\ }|^{{\ }^{\!\!\!\!-\!\!-\!\!-\!\!-\!\!-\!\!-\!\!-\!\!}}_{\ \
\sqcap}|\ \ \ |\cr
 1\ \ 2\ 3\, 4\ 5\ \ 6\end{matrix}\right\}$$

As for the corresponding element of $NC(5,1)$, this is pictured as follows:
$$p_{51}=\left\{\begin{matrix}
1\ \ 2\  3\, 4\ 5\cr
|_{{\ }_{\!\!\!\!-\!\!-\!\!-\!\!-\!\!}}^
{\ }|_{{\ }_{\!\!\!\!-\!\!-\!\!-\!\!-\!\!-\!\!-\!\!-\!\!}}^{\ \ \sqcup}|\cr
{\ }_|\cr
\ 1
\end{matrix}\right\}$$

We fix now a number $n\in\mathbb N$. All indices will vary in the set $\{1,\ldots,n\}$.

\begin{definition}
Asociated to any partition $p\in NC(k,l)$ and any multi-indices $i=(i_1,\ldots,i_k)$ and $j=(j_1,\ldots,j_l)$ is a number $p(i,j)\in\{0,1\}$, as follows:
\begin{enumerate}
\item We put the indices of $i,j$ on the points of $p$, in the obvious way.

\item If all the strings of $p$ join equal indices, we set $p(i,j)=1$.

\item If some strings of $p$ join different indices, we set $p(i,j)=0$.
\end{enumerate}
\end{definition}

Here is a series of basic examples, with the partitions represented by the corresponding pictures, drawn according to the above conventions:
$$\left\{\begin{matrix}1\cr|\cr1\end{matrix}\right\}(a,b)
=\left\{\begin{matrix}\ \cr\sqcap\cr 1\,2\end{matrix}\right\}(,ab)
=\left\{\begin{matrix}1\,2 \cr\sqcup\cr\ \end{matrix}\right\}(ab,)
=\delta_{ab}$$

In this equality the $\delta$ symbol on the right is a usual Kronecker symbol.

\begin{definition}
Associated to any partition $p\in NC(k,l)$ is the linear map
$$T_p(e_{i_1}\otimes\ldots\otimes e_{i_k})=\sum_{j_1\ldots j_l}p(i,j)\,e_{j_1}\otimes\ldots\otimes e_{j_l}$$
where $e_1,\ldots,e_n$ is the standard basis of $\mathbb C^n$.
\end{definition}

Here are a few examples, of interest for the considerations to follow: 
$$T\underline{\ }\left\{\big\vert\ \big\vert\right\}(e_a\otimes e_b)=e_a\otimes e_b$$
$$T\underline{\ }\left\{\big\vert\!\!\!-\!\!\!\big\vert\right\}(e_a\otimes e_b)=\delta_{ab}\,e_a\otimes e_a$$
$$T\underline{\ }\left\{\!\!\begin{matrix}{\ }^{|\ |}\cr {\ }_{|\ |}\end{matrix}\right\}(e_a\otimes e_b)=\sum_{cd}e_c\otimes e_d$$
$$T\underline{\ }\left\{\!\!\begin{matrix}\ \sqcup\cr{\ }_{|\ |}\end{matrix}\right\}(e_a\otimes e_b)=\delta_{ab}\sum_{cd}e_c\otimes e_d$$
$$T\underline{\ }\left\{\begin{matrix}\sqcup\cr\sqcap\end{matrix}\right\}(e_a\otimes e_b)=\delta_{ab}\sum_ce_c\otimes e_c$$

We introduce now a number of algebraic operations on partitions.

\begin{definition}
The tensor product, composition and involution of partitions are obtained by horizontal and vertical concatenation and upside-down turning
\begin{eqnarray*}
p\otimes q&=&\{\mathcal P\mathcal Q\}\\
pq&=&\left\{\begin{matrix}\mathcal Q\cr \mathcal P\end{matrix}\right\}-\{{\rm closed\ blocks}\}\\
p^*&=&\{\mathcal P^\curvearrowright\}
\end{eqnarray*}
where $p=\{\mathcal P\}$ and $q=\{\mathcal Q\}$ are the pictorial representations of $p,q$.
\end{definition}

The above three operations can be of course defined by certain explicit algebraic formulae, by using the formalism of partitions, but we prefer to use their pictorial interpretation, which is less heavier, and far more suggestive.

As an example, consider two partitions $p\in NC(k,l)$ and $q\in NC(k',l')$. In order to define their tensor product, we use the following identification:
$$\{1,\ldots,k,1,\ldots,l\}\sqcup\{1,\ldots,k',1,\ldots,l'\}\simeq \{1,\ldots,k+k',1,\ldots,l+l'\}$$

Here the elements of the first set on the left are identified with the corresponding elements on the set on the right, and the elements of the second set on the left are identified with the missing elements at right, in the following way:
$$\{1,\ldots,k',1,\ldots,l'\}\simeq \{k+1,\ldots,k+k',l+1,\ldots,l+l'\}$$

Now with the above identification, the disjoint union $p\sqcup q$ is a partition of the union of the two sets on the left, hence can be regarded as a partition of the big set on the right. We denote this latter partition by $p\otimes q\in NC(k+k',l+l')$.

Observe that in pictorial terms, this partition $p\otimes q$ is simply obtained by ``horizontal concatenation'', as stated in Definition 5.4.

The composition is similarly defined by ``vertical concatenation'' of the pictures. Observed that it is only partially defined: the number of upper points of $p$ must be equal to the number of lower points of $q$. Moreover, when identifying the upper points of $p$ with the lower points of $q$, ``closed blocks'' might appear, i.e. strings which are not connected to any of the new upper and lower points. These blocks are simply discarded from the concatenated picture.

We are now in position of developing the method explained in the beginning of this section. Let $U$ be the fundamental corepresentation of $A_s(n)$.

\begin{theorem}
We have the equality
$$Hom(U^{\otimes k},U^{\otimes l})={\rm span}\{T_p|p\in NC(k,l)\}$$ 
and if $n\geq 4$, the maps on the right are linearly independent.
\end{theorem}

\begin{proof}
This result is known since \cite{ba2}, a simplified proof being as follows. First, it is routine to check that we have the following formulae, with $b(p,q)\in\mathbb N$:
\begin{eqnarray*}
T_{p\otimes q}&=&T_p\otimes T_q\\
T_{pq}&=&n^{-b(p,q)}T_pT_q\\
T_{p^*}&=&T_p^*\\
T_{|}&=&id
\end{eqnarray*}

This shows that the spaces on the right form a tensor category in the sense of Woronowicz \cite{wo2}. Moreover, since $T_\sqcap$ implements the canonical ``duality'' map, this category has duals, so by \cite{wo2} it gives rise to a certain Hopf algebra $(A,u)$.

Now since the one-block partitions $1_k\in NC(k)$ produce via the operations in Definition 5.4 all the noncrossing partitions, our tensor category is generated by the maps $T_{1_k}$. This means that the algebra $A$ is presented by the relations $T_{1_k}\in Hom(1,u^{\otimes k})$, and a routine computation shows that these conditions are equivalent to the fact that $u$ is magic. Thus we have $A=A_s(n)$.

The last assertion, proved as well in \cite{ba2}, follows by using a number of standard tricks. First, by Frobenius duality, the validity of the statement depends only on the value of $k+l$. Moreover, once again by a standard representation theory argument, coming from $1\in U$, we can assume that $k+l$ is even. Thus it is enough to do the check in the case $k=l$. But here the vector spaces in the statement are actually algebras, and the result follows by using a suitable positive trace.
\end{proof}

The above result is not suprising, because for $C(S_n)$ the corresponding spaces of intertwiners are given by the same formula, but with all the partitions instead of just the noncrossing ones. Thus we are in  tune with the general principle ``the passage from classical to free is obtained by restricting attention to the noncrossing partitions'', which goes back to Speicher's paper \cite{spe}.

\section{Tannakian duality}

We are now in position of investigating the spaces of intertwiners between the various tensor products of basic corepresentations of $A_h^s(n)$. 

\begin{definition}
We make the assumption $n\geq 4$.
\end{definition}

This assumption, to be kept until the end of the paper, guarantees that the linear maps $T_p$ in Theorem 5.5 are linearly independent. We will use this fact in order to identify the partitions with the corresponding linear maps.

In the cases $n=2,3$, not to be investigated in what follows, $A_h^s(n)$ collapses to a quite simple algebra, as shown by Corollary 3.5.

\begin{definition}
For $i_1,\ldots,i_k\in\mathbb Z$ we use the notation
$$u_{i_1\ldots i_k}=u_{i_1}\otimes\ldots\otimes u_{i_k}$$
where $\{u_i|i\in\mathbb Z\}$ are the corepresentations in Theorem 4.3.
\end{definition}

Observe that in the particular case $i_1,\ldots,i_k\in\{\pm 1\}$, we obtain in this way all the possible tensor products between $u=u_1$ and $\bar{u}=u_{-1}$, known by the general results in \cite{wo1} to contain any irreducible corepresentation of $A_h^s(n)$.

\begin{theorem}
We have the equality
$$Hom(u_{i_1\ldots i_k},u_{j_1\ldots j_l})={\rm span}\{T_p|p\in NC_s(i_1\ldots i_k,j_1\ldots j_l)\}$$
where the set on the right consists of elements of $NC(k,l)$ having the property that in each block, the sum of $i$ indices equals the sum of $j$ indices, modulo $s$.
\end{theorem}

\begin{proof}
The idea will be to expand, suitably modify, and unify the proof of the following key particular cases:
\begin{enumerate}
\item Theorem 5.5 from the previous section, which gives the result at $s=1$, for any choice of indices $i_1,\ldots,i_k,j_1,\ldots,j_l\in\mathbb Z$.

\item The main technical result in \cite{bb+}, page 37 on top there, which gives the result for any $s\in\mathbb N$, for indices of type $i_1,\ldots,i_k,j_1,\ldots,j_l\in\{\pm 1\}$. 
\end{enumerate}

{\bf Step 1.} Our first claim is that, in order to prove $\supset$, we may restrict attention to the case $k=0$. Indeed, it is known that for any two corepresentations $v,w$ we have a Frobenius duality isomorphism, as follows:
$$Hom(v,w)\simeq Hom(1,v\otimes\bar{w})$$

In the case $v=u_{i_1\ldots i_k}$ and $w=u_{j_1\ldots j_l}$, we can use the formulae in Theorem 4.3 in order to compute $v\otimes\bar{w}$, and the Frobenius isomorphism reads:
$$Hom(u_{i_1\ldots i_k},u_{j_1\ldots j_l})\simeq Hom(1,u_{i_1\ldots i_k(-j_l)\ldots(-j_1)})$$

On the other hand, we have the following canonical identification:
$$NC(k,l)\simeq NC(0,k+l)$$

Now it follows from definitions and from Theorem 5.5 that at $s=1$ these two isomorphisms are compatible, in the obvious sense. Together with the functoriality considerations regarding the canonical map $A_h^s(n)\to A_s(n)$, explained in the beginning of the previous section, this justifies our claim.

{\bf Step 2.} Our second claim is that, in order to prove $\supset$ in the case $k=0$, we may restrict attention to the one-block partitions. Indeed, this follows once again from a standard trick. Consider the following disjoint union:
$$NC_s=\bigcup_{k=0}^\infty\bigcup_{i_1\ldots i_k} NC_s(0,i_1\ldots i_k)$$

This is a set of labeled partitions, having the following properties:
\begin{enumerate}
\item Each $p\in NC_s$ is noncrossing.

\item For $p\in NC_s$, any block of $p$ is in $NC_s$. 
\end{enumerate}

It is well-known that under these assumptions, the global algebraic properties of $NC_s$ can be checked on blocks, and this justifies our claim.

{\bf Step 3.} We finish the proof of $\supset$. According to the above considerations, we just have to prove that the vector associated to the one-block partition in $NC(l)$ is fixed by $u_{j_1\ldots j_l}$, for any choice of $j_1,\ldots ,j_l$ satisfying:
$$s|j_1+\ldots+j_l$$

Consider the standard generators $e_{ab}\in M_n(\mathbb C)$, acting on the basis vectors by $e_{ab}(e_c)=\delta_{bc}e_a$. The corepresentation $u_{j_1\ldots j_l}$ is given by the following formula:
\begin{eqnarray*}
u_{j_1\ldots j_l}
&=&u_{j_1}\otimes\ldots\otimes u_{j_l}\\
&=&\left(u_{a_1b_1}^{j_1}\right)\otimes\ldots\otimes\left(u_{a_lb_l}^{j_l}\right)\\
&=&\sum_{a_1\ldots a_l}\sum_{b_1\ldots b_l}u_{a_1b_1}^{j_1}\ldots u_{a_lb_l}^{j_l}\otimes e_{a_1b_1}\otimes\ldots\otimes e_{a_lb_l}
\end{eqnarray*}

As for the vector associated to the one-block partition, this is:
$$\xi_l=\sum_be_b^{\otimes l}$$

By using several times the relations in Lemma 4.2, we get as claimed: 
\begin{eqnarray*}
u_{j_1\ldots j_l}(1\otimes\xi_l)
&=&\sum_{a_1\ldots a_l}\sum_bu_{a_1b}^{j_1}\ldots u_{a_lb}^{j_l}\otimes e_{a_1}\otimes\ldots\otimes e_{a_l}\\
&=&\sum_{ab}u_{ab}^{j_1+\ldots+j_l}\otimes e_a^{\otimes l}\\
&=&\sum_{ab}p_{ab}\otimes e_a^{\otimes l}\\
&=&\sum_a1\otimes e_a^{\otimes l}\\
&=&1\otimes\xi_l
\end{eqnarray*}

{\bf Step 4.} We begin the proof of $\subset$. The first remark, which can be justified as in the proof of Theorem 5.5, is that the spaces on the right in the statement form a tensor category with duals in the sense of Woronowicz \cite{wo2}. Thus by Tannakian duality they correspond to a certain Hopf algebra $A$. 

This algebra is by definition the maximal model for the tensor category. In other words, it comes with a family of corepresentations $\{v_i\}$, such that:
$$Hom(v_{i_1\ldots i_k},v_{j_1\ldots j_l})={\rm span}\{T_p|p\in NC_s(i_1\ldots i_k,j_1\ldots j_l)\}$$

Here, and in what follows, we use the notation $v_{i_1\ldots i_k}=v_{i_1}\otimes\ldots\otimes v_{i_k}$.

On the other hand, the inclusion $\supset$ that we just proved shows that $A_h^s(n)$ is a model for the tensor category. Thus by \cite{wo2} we have a surjective arrow $A\to A_h^s$, mapping $v_i\to u_i$ for any $i$. We have to prove that this is an isomorphism.

{\bf Step 5.} We finish the proof of $\subset$. This can be done in a straightforward way, by suitably adapting the proof of Theorem 5.5. In what follows we present a shorter argument, based on some previous work in \cite{bb+}. The main technical result there (page 37, on top) is that the equality in the statement holds in the case where both $u_{i_1\ldots i_k}$ and $u_{j_1\ldots j_l}$ are tensor products between $u$ and $\bar{u}$. 

With the present notations we have $u=u_1$ and $\bar{u}=u_{-1}$, so what we know from \cite{bb+} is that the result holds for any choice of indices $i_r,j_r\in\{\pm 1\}$. 

Now by using the definition of $A$, we get that for such indices we have:
$$Hom(u_{i_1\ldots i_k},u_{j_1\ldots j_l})
=Hom(v_{i_1\ldots i_k},v_{j_1\ldots j_l})$$

In other words, the map $A\to A_h^s(n)$ induces isomorphisms at the level of intertwining spaces between the various tensor products between $v$ and $\bar{v}$. It is well-known that such a map must be an isomorphism, and this finishes the proof.
\end{proof}

As an illustration for the above result, we present below two corollaries, both of them with very detailed proofs. 

First is a key statement about the basic corepresentations $u_i$. As usual, we use indices modulo $s$, with the convention $u_{ij}^0=p_{ij}$.

\begin{corollary}
The basic corepresentations $u_0,\ldots,u_{s-1}$ are as follows:
\begin{enumerate}
\item $u_1,\ldots,u_{s-1}$ are irreducible.

\item $u_0=1+r_0$, with $r_0$ irreducible.

\item $r_0,u_1,\ldots,u_{s-1}$ are distinct.
\end{enumerate}
\end{corollary}

\begin{proof}
We apply Theorem 6.3 with $k=l=1$ and $i_1=i,j_1=j$. This gives:
$$\dim (Hom(u_i,u_j))=\# NC_s(i,j)$$

We have two candidates for the elements of $NC_s(i,j)$, namely the two partitions in $NC(1,1)$. So, consider these two partitions, with the points labeled by $i,j$:
$$p=\left\{\begin{matrix}i\cr\Big\vert\cr j\end{matrix}\right\}\quad\quad\quad
q=\left\{\begin{matrix}i\cr|\cr\cr |\cr j\end{matrix}\right\}$$

We have to check for each of these partitions if the sum of $i$ indices equals or not the sum of $j$ indices, modulo $s$, in each block. The answer is as follows:
\begin{eqnarray*}
p\in NC_s(i,j)&\iff&i=j\\
q\in NC_s(i,j)&\iff&i=j=0
\end{eqnarray*}

By collecting together these two answers, we get:
$$\# NC_s(i,j)=
\begin{cases}
0&{\rm if\ }i\neq j\\
1&{\rm if\ }i=j\neq 0\\
2&{\rm if\ }i=j=0
\end{cases}$$

This gives all the results. Indeed, (1) follows from the second equality, (2) follows from the third equality and from the fact that we have $1\in u_s$ (this is because $u_s=p$ is magic), and (3) follows from the first equality.

As a last remark, the ingredient $1\in u_0$ can be deduced as well from Theorem 6.3. Indeed, with $k=0,l=1$ and $j_1=0$ we get $\dim(Hom(1,u_0))=1$.
\end{proof}

The second corollary is a key statement, to be used in what follows for the computation of the fusion rules.

It is convenient at this point to switch back to the old notation for the tensor products between the basic corepresentations, by ignoring Definition 6.2, which was temporary. Also, we use the notation $\#(1\in v)=\dim(Hom(1,v))$.

\begin{corollary}
We have the formula
$$\#(1\in u_{i_1}\otimes\ldots\otimes u_{i_k})=\# NC_s(i_1\ldots i_k)$$
where the set on the right consists of noncrossing partitions of $\{1,\ldots,k\}$ having the property that the sum of indices in each block is a multiple of $s$.
\end{corollary}

\begin{proof}
This is clear from Theorem 6.3.
\end{proof}

\section{The main result}

It is known from Woronowicz's analogue of Peter-Weyl theory in \cite{wo1} that each corepresentation decomposes as a direct sum of irreducible corepresentations. 

In particular any tensor product of irreducible corepresentations decomposes as a direct sum of irreducible corepresentations.

The formulae describing these decompositions are called fusion rules.

\begin{definition}
The fusion semiring $(R^+,-,+,\otimes)$ is defined as follows:
\begin{enumerate}
\item $R^+$ is the set of equivalence classes of corepresentations.

\item $-,+,\otimes$ are the usual involution, sum and tensor product.
\end{enumerate}
\end{definition}

It follows from the Peter-Weyl type results that $(R^+,+)$ is the free additive monoid on the set of irreducible corepresentations. Thus the fusion semiring $(R^+,-,+,\otimes)$ encodes the collection of fusion rules.

\begin{definition}
Let $F=<\mathbb Z_s>$ be the monoid formed by the words over $\mathbb Z_s$. We endow $F$ with the following operations:
\begin{enumerate}
\item Involution: $(i_1\ldots i_k)^-=(-i_k)\ldots(-i_1)$.

\item Fusion: $(i_1\ldots i_k)\cdot (j_1\ldots j_l)=i_1\ldots i_{k-1}(i_k+j_1)j_2\ldots j_l$.
\end{enumerate}
\end{definition}

Note that $v\cdot w$ is not defined when $v$ or $w$ is the empty word. We make the convention that the corresponding terms disappear from the fusion rules below.

We are now in position of stating the main result in this paper. We recall from Corollary 6.4 that the basic corepresentations $u_1,\ldots,u_s$ are all irreducible, except for the last one, which is of the form $1+r_s$, with $r_s$ irreducible.

\begin{theorem}
The irreducible corepresentations of $A_h^s(n)$ can be labeled $r_x$ with $x\in F$, such that the involution and fusion rules are $\bar{r}_x=r_{\bar{x}}$ and
$$r_x\otimes r_y=\sum_{x=vz,y=\bar{z}w}r_{vw}+r_{v\cdot w}$$
and such that we have $r_i=u_i-\delta_{i0}1$ for any $i\in\mathbb Z_s$.
\end{theorem}

\begin{proof}
We use the results in the previous section, and a standard method from \cite{ba1}, \cite{ba2}. Consider the set of irreducible corepresentations of $A_h^s(n)$, its fusion semiring, and its fusion ring, and denote them as follows:
$$R_{irr}\subset R^+\subset R$$

Observe that by Corollary 6.4, we have $r_i\in R_{irr}$ for any $i$.

{\bf Step 1.} We first construct an abstract fusion semiring, with fusion rules as in the statement. Consider indeed the monoid $A=\{a_x|x\in F\}$, with multiplication $a_xa_y=a_{xy}$. We denote by $\mathbb NA$ the set of linear combinations of elements in $A$, with coefficients in $\mathbb N$, and we endow it with fusion rules as in the statement:
$$a_x\otimes a_y=\sum_{x=vz,y=\bar{z}w}a_{vw}+a_{v\cdot w}$$

With these notations, $(\mathbb NA,+,\otimes)$ is a semiring. We will use as well the set $\mathbb ZA$, formed by the linear combinations of elements of $A$, with coefficients in $\mathbb Z$. The above tensor product operation extends to $\mathbb ZA$, and $(\mathbb ZA,+,\otimes)$ is a ring.

{\bf Step 2.} We claim that the fusion rules on $\mathbb ZA$ can be uniquely described by conversion formulae as follows:
$$a_{i_1}\otimes\ldots\otimes a_{i_k}=\sum_l\sum_{j_1\ldots j_l}C_{i_1\ldots i_k}^{j_1\ldots j_l}a_{j_1\ldots j_l}$$
$$a_{i_1\ldots i_k}=\sum_l\sum_{j_1\ldots j_l}D_{i_1\ldots i_k}^{j_1\ldots j_l}a_{j_1}\otimes\ldots\otimes a_{j_l}$$

Here the $C$ coefficients are certain positive integers, and the $D$ coefficients are certain integers. The existence and uniqueness of such decompositions follow indeed from the definition of the tensor product operation, and by induction over $k$ for the $D$ coefficients.

{\bf Step 3.} We claim that there is a unique morphism of rings $\Phi:\mathbb ZA\to R$, such that $\Phi(a_i)=r_i$ for any $i$. Indeed, consider the following elements of $R$:
$$r_{i_1\ldots i_k}=\sum_l\sum_{j_1\ldots j_l}D_{i_1\ldots i_k}^{j_1\ldots j_l}r_{j_1}\otimes\ldots\otimes r_{j_l}$$

In case we have a morphism as claimed, by linearity and multiplicativity we must have $\Phi(a_x)=r_x$ for any $x\in F$. Thus our morphism is uniquely determined on  $A$, so by linearity it is uniquely determined on $\mathbb ZA$.

In order to prove the existence, we can set $\Phi(a_x)=r_x$ for any $x\in F$, then extend $\Phi$ by linearity to the whole $\mathbb ZA$. Since $\Phi$ commutes with the above conversion formulae, which describe the fusion rules, it is indeed a morphism. 

{\bf Step 4.} We claim that $\Phi$ commutes with the linear forms $x\to\#(1\in x)$. Indeed, by linearity we just have to check the following equality:
$$\#(1\in a_{i_1}\otimes\ldots\otimes a_{i_k})=\#(1\in r_{i_1}\otimes\ldots\otimes r_{i_k})$$

Now remember that the elements $r_i$ are defined as $r_i=u_i-\delta_{i0}1$. So, consider the elements $c_i=a_i+\delta_{i0}1$. Since the operations $r_i\to u_i$ and $a_i\to c_i$ are of the same nature, by linearity the above formula is equivalent to:
$$\#(1\in c_{i_1}\otimes\ldots\otimes c_{i_k})=\#(1\in u_{i_1}\otimes\ldots\otimes u_{i_k})$$

Now by using Corollary 6.5, what we have to prove is:
$$\#(1\in c_{i_1}\otimes\ldots\otimes c_{i_k})=\#NC_s(i_1\ldots i_k)$$

In order to prove this formula, consider the product on the left:
$$P=(a_{i_1}+\delta_{i_10}1)\otimes(a_{i_2}+\delta_{i_20}1)\otimes\ldots\otimes (a_{i_k}+\delta_{i_k0}1)$$

This quantity can be computed by using the fusion rules on $A$. An induction on $k$ shows that the final components of type $a_x$ will come from the different ways of grouping and summing the consecutive terms of the sequence $(i_1,\ldots,i_k)$, and simultaneously removing some of the sums which vanish modulo $s$, so as to obtain the sequence $x$. This can be encoded by families of noncrossing partitions, and in particular the 1 components will come from the partitions in $NC_s(i_1\ldots i_k)$. Thus we have $\#(1\in P)=\# NC_s(i_1\ldots i_k)$, as claimed.

{\bf Step 5.} We claim that $\Phi$ is injective. Indeed, this follows from the result in the previous step, by using a standard positivity argument:
\begin{eqnarray*}
\Phi(\alpha)=0
&\implies&\Phi(\alpha\alpha^*)=0\\
&\implies&\#(1\in \Phi(\alpha\alpha^*))=0\\
&\implies&\#(1\in \alpha\alpha^*)=0\\
&\implies&\alpha=0
\end{eqnarray*}

Here $\alpha$ is arbitrary in the domain of $\Phi$, we use the notation $a_x^*=a_{\bar{x}}$, where $x\to\bar{x}$ is the involution in Definition 7.2, and $a\to\#(1,a)$ is the unique linear extension of the operation consisting of counting the number of 1's.  Observe that this latter linear form is indeed positive definite, according to the identity $\#(1,a_xa_y^*)=\delta_{xy}$, which is clear from the definition of the product of $\mathbb ZA$.

{\bf Step 6.} We claim that we have $\Phi(A)\subset R_{irr}$. This is the same as saying that $r_x\in R_{irr}$ for any $x\in F$, and we will prove it by recurrence on the length of $x$. 

For the words of length $1$ the assertion is true because of Corollary 6.4, as pointed out in the beginning of the proof.

So, assume that the assertion is true for all the words of length $<k$, and consider an arbitrary length $k$ word, $x=i_1\ldots i_k$. We have:
$$a_{i_1}\otimes a_{i_2\ldots i_k}=a_x+a_{i_1+i_2,i_3\ldots i_k}+\delta_{i_1+i_2,0}a_{i_3\ldots i_k}$$

By applying $\Phi$ to this decomposition, we get:
$$r_{i_1}\otimes r_{i_2\ldots i_k}=r_x+r_{i_1+i_2,i_3\ldots i_k}+\delta_{i_1+i_2,0}r_{i_3\ldots i_k}$$

For $u$ irreducible, we use the notation $\#(u\in v)=\dim(Hom(u,v))$. We have the following computation, which is valid for $y=i_1+i_2,i_3\ldots i_k$, as well as for $y=i_3\ldots i_k$ in the case $i_1+i_2=0$:
\begin{eqnarray*}
\#(r_y\in r_{i_1}\otimes r_{i_2\ldots i_k})
&=&\#(1,r_{\bar{y}}\otimes r_{i_1}\otimes r_{i_2\ldots i_k})\\
&=&\#(1,a_{\bar{y}}\otimes a_{i_1}\otimes a_{i_2\ldots i_k})\\
&=&\#(a_y\in a_{i_1}\otimes a_{i_2\ldots i_k})\\
&=&1  
\end{eqnarray*}

Moreover, we know from the previous step that we have $r_{i_1+i_2,i_3\ldots i_k}\neq r_{i_3\ldots i_k}$, so we conclude that the following formula defines an element of $R^+$:
$$\alpha=r_{i_1}\otimes r_{i_2\ldots i_k}-r_{i_1+i_2,i_3\ldots i_k}-\delta_{i_1+i_2,0}r_{i_3\ldots i_k}$$

On the other hand, we have $\alpha=r_x$, so we conclude that we have $r_x\in R^+$. Finally, the irreducibility of $r_x$ follows from the following computation:
\begin{eqnarray*}
\#(1\in r_x\otimes\bar{r}_x)
&=&\#(1\in r_x\otimes r_{\bar{x}})\\
&=&\#(1\in a_x\otimes a_{\bar{x}})\\
&=&\#(1\in a_x\otimes\bar{a}_x)\\
&=&1
\end{eqnarray*}

{\bf Step 7.} Summarizing, we have constructed an injective ring morphism $\Phi:\mathbb ZA\to R$, having the property $\Phi(A)\subset R_{irr}$. The remaining fact to be proved, namely that we have $\Phi(A)=R_{irr}$, is clear from the general results in \cite{wo1}. Indeed, since each element of $\mathbb NA$ is a sum of elements in $A$, by applying $\Phi$ we get that each element in $\Phi(\mathbb NA)$ is a sum of irreducible corepresentations in $\Phi(A)$. But since $\Phi(\mathbb NA)$ contains all the tensor powers between the fundamental corepresentation and its conjugate, we get by \cite{wo1} that we have $\Phi(A)=R_{irr}$, and we are done.
\end{proof}

As an illustration for the above result, we would like to work out the case $s=1$. We get a new proof for the following result from \cite{ba2}.

\begin{corollary}
The fusion rules for $A_s(n)$ are the same as the Clebsch-Gordan rules for the irreducible representations of $SO_3$.
\end{corollary}

\begin{proof}
We first examine Definition 7.2, in the case $s=1$. We have a canonical isomorphism $F\simeq\mathbb N$, mapping the words to their lengths. With this identification, the involution is $\bar{k}=k$, and the fusion is $k\cdot l=k+l-1$. 

We apply now Theorem 7.3. The conclusion is that the irreducible corepresentations can be labeled $\{r_k|k\in\mathbb N\}$, and that we have:
\begin{enumerate}
\item $\bar{r}_k=r_k$.

\item $r_k\otimes r_l=r_{k+l}+r_{k+l-1}+\ldots+r_{|k-l|}$.

\item $r_1=u-1$.
\end{enumerate}

In other words, we have reached to the Clebsch-Gordan rules for $SO_3$, with $u-1$ corresponding to the 3-dimensional representation of $SO_3$.
\end{proof}

\section{Alternative formulation}

In this section we present an alternative formulation of Theorem 7.3, that we will once again illustrate with a complete computation at $s=1$.

We begin with a slight modification of Theorem 7.3.

\begin{theorem}
Consider the free monoid $A=<a_i|i\in\mathbb Z_s>$ with the involution $a_i^*=a_{-i}$, and define inductively the following fusion rules on it:
$$pa_i\otimes a_jq=pa_ia_jq+pa_{i+j}q+\delta_{i+j,0}p\otimes q$$

Then the irreducible corepresentations of $A_h^s(n)$ can be indexed by the elements of $A$, and the fusion rules and involution are the above ones.
\end{theorem}

\begin{proof}
We claim that this follows from Theorem 7.3, by performing the following relabeling of the irreducible corepresentations:
$$r_{i_1\ldots i_k}\to a_{i_1}\ldots a_{i_k}$$

Indeed, with the notations in Theorem 7.3 we have the following computation, valid for any two elements $i,j\in\mathbb Z_s$ and any two words $x,y\in F$:
\begin{eqnarray*}
r_{xi}\otimes r_{jy}
&=&\sum_{xi=vz,jy=\bar{z}w}r_{vw}+r_{v\cdot w}\\
&=&r_{xijy}+r_{x,i+j,y}+\delta_{i+j,0}\sum_{x=vz,y=\bar{z}w}r_{vw}+r_{v\cdot w}\\
&=&r_{xijy}+r_{x,i+j,y}+\delta_{i+j,0}r_x\otimes r_y
\end{eqnarray*}

With the above relabeling $r_{i_1\ldots i_k}\to a_{i_1}\ldots a_{i_k}$, this gives the formula in the statement (with $r_x\to p$ and $r_y\to q$), and we are done.
\end{proof}

Our alternative reformulation of Theorem 7.3 is based on the idea of embedding $R^+$ into a bigger fusion semiring. In this bigger semiring the fusion rules will appear to be actually simpler, due to a Fourier transform type situation.

Given a fusion monoid $M$ and an element $b\in M$, we denote by $<b>$ the fusion monoid generated by $b$. In other words, $<b>$ is the smallest subset of $M$ containing $b$, and which is stable by composition, involution and fusion rules. 

\begin{theorem}
Consider the monoid $M=<a,z|z^s=1>$ with the involution $a^*=a,z^*=z^{-1}$, and define inductively the following fusion rules on it:
$$vaz^i\otimes z^jaw=vaz^{i+j}aw+\delta_{s|i+j}v\otimes w$$

Then the irreducible corepresentations of $A_h^s(n)$ can be indexed by the elements of $N=<aza>$, and the fusion rules and involution are the above ones.
\end{theorem}

\begin{proof}
It is routine to check that the elements $az^ia$ with $i=1,\ldots,s$ are free inside $M$. In other words, the submonoid $N'=<az^ia>$ is free on $s$ generators, so it can be identified with the free monoid $A$ in Theorem 8.1, via $a_i=az^ia$.

We have $(az^ia)^*=az^{-i}a$, so this identification is involution-preserving.

Consider now two arbitrary elements $p,q\in N'$. By using twice the formula in the statement, we get the formula in Theorem 8.1:
\begin{eqnarray*}
pa_i\otimes a_jq
&=&paz^ia\otimes az^jaq\\
&=&paz^iaaz^jaq+paz^i\otimes z^jaq\\
&=&paz^iaaz^jaq+paz^{i+j}aq+\delta_{i+j,0}p\otimes q\\
&=&pa_ia_jq+pa_{i+j}q+\delta_{i+j,0}p\otimes q
\end{eqnarray*}

Thus our identification $N'\simeq A$ is fusion rule-preserving.

In order to conclude, it remains to prove that the inclusion $N\subset N'$ is actually an equality. But this follows from the fact that $A$ is generated as a fusion monoid by $a_1$. Indeed, by using the identification $N'\simeq A$ this shows that $N'$ is generated as a fusion monoid by $aza$, and we are done. 
\end{proof}

As in illustration for the above result, we work out the case $s=1$, leading to a supplementary proof for the main result in \cite{ba2}.

\begin{corollary}
At $s=1$ we have the Clebsch-Gordan rules for $SO_3$.
\end{corollary}

\begin{proof}
Indeed, at $s=1$ we have a canonical isomorphism $M\simeq\mathbb N$, and with this identification, the fusion rules are given inductively by:
$$k\otimes l=(k+l)+(k-1)\otimes (l-1)$$

This gives the following explicit formula, which is nothing but the Clebsch-Gordan formula for the fusion rules of irreducible representations of $SU_2$:
$$k\otimes l=(k+l)+(k+l-2)+\ldots+|k-l|$$

As for the submonoid $N\subset M$, this corresponds via the above identifications to the submonoid $2\mathbb N\subset \mathbb N$. Thus when performing a division by 2 we get an isomorphism $N\simeq\mathbb N$, and we have the following fusion rules on $N$:
$$k\otimes l=(k+l)+(k+l-1)+\ldots+|k-l|$$

We recognize here the Clebsch-Gordan formula for $SO_3$.
\end{proof}

Let us also work out the case $s=\infty$.

\begin{corollary}
Consider the monoid $M=<a,z,z^{-1}>$ with the involution $a^*=a,z^*=z^{-1}$, and define inductively the following fusion rules on it:
$$vaz^i\otimes z^jaw=vaz^{i+j}aw+\delta_{i+j,0}v\otimes w$$

Then the irreducible corepresentations of $A_h^\infty(n)$ can be indexed by the elements of $N=<aza>$, and the fusion rules and involution are the above ones.
\end{corollary}

\begin{proof}
This is a reformulation of Theorem 8.2 in the case $s=\infty$, by using the various conventions and notations specific to this case.
\end{proof}

\section{Dimension formula}

In this section we compute the dimension of the irreducible corepresentations of $A_h^s(n)$. Besides of being of independent theoretical interest, this computation can be regarded as being first ingredient towards a fine study of the growth invariants of $A_h^s(n)$, in the spirit of \cite{ve1}, \cite{ve2}, \cite{bve}.

The main result will be best expressed in terms of the alternative formalism from the previous section. So, consider the free monoids $N\subset M$ in Theorem 8.2, with $N=<aza>$ labelling the irreducible corepresentations of $A_h(n)$.

We have seen in the proof of Corollary 8.3 that, at $s=1$, the inclusion $N\subset M$ ultimately comes from double cover map $SU_2\to SO_3$. So, our first task will be to introduce a certain sequence of numbers $d_k$, which appear as ``versions, with $2$ replaced by $\sqrt{n}$'' of the dimensions of the irreducible representations of $SU_2$.

\begin{definition}
Associated to any $n\geq 4$ is the sequence of numbers given by 
$$d_{k+1}+d_{k-1}=\sqrt{n}d_k$$
with the initial values $d_0=1$ and $d_1=\sqrt{n}$.
\end{definition}

These numbers can be of course computed in terms of the roots of the polynomial $X^2-\sqrt{n}X+1=0$, but the explicit formula is not very enlightening. Instead, let us just indicate the first few values of these numbers:
\begin{eqnarray*}
d_0&=&1\\
d_1&=&\sqrt{n}\\
d_2&=&n-1\\
d_3&=&(n-2)\sqrt{n}\\
d_4&=&n^2-3n+1
\end{eqnarray*}

Observe that at $n=4$ the recurrence relation is $d_{k+1}+d_{k-1}=2d_k$, with initial values $d_0=1$ and $d_1=2$, so we have $d_k={k+1}$ for any $k$. These numbers are indeed the dimensions of the irreducible representations of $SU_2$.

In the general case, the numbers $d_k$ appear as dimensions of the irreducible corepresentations of the Wang algebra $A_o(\sqrt{n})$. Without getting into details here (see \cite{wa1}), let us just record the following useful fact.

\begin{proposition}
The sequence of numbers $d_k$ can be defined alternatively by the Clebsch-Gordan type formula
$$d_kd_l=d_{k+l}+d_{k-1}d_{l-1}$$
with the initial values $d_0=1$ and $d_1=\sqrt{n}$.
\end{proposition}

\begin{proof}
Let $M\simeq\mathbb N$ be endowed with the Clebsch-Gordan rules, as in the proof of Corollary 8.3. It follows from definitions that associated to any $\lambda\in\mathbb R$ is a unique morphism of fusion semirings $d:\mathbb NM\to (\mathbb N,+,\cdot)$, having the property $d_1=\lambda$. According to the Clebsch-Gordan rules, we must have:
$$d_kd_l=d_{k+l}+d_{k+l-2}+\ldots+d_{|k-l|}$$

This shows that the numbers $d_k$ satisfy the recurrence formula in the statement, so if we set $\lambda=\sqrt{n}$ we get indeed the numbers in the statement. 

On the other hand, with $l=1$, the above formula reads:
$$d_k\sqrt{n}=d_{k+1}+d_{k-1}$$

Together with $d_0=1$ and $d_1=\sqrt{n}$, this shows that our sequence of numbers $d_k$ coincides with the one in Definition 9.1, and we are done.
\end{proof}

\begin{theorem}
The dimensions of the irreducible corepresentations of $A_h^s(n)$, as labelled by the monoid $N=<aza>$ in Theorem 8.2, are given by
$$\dim(a^{i_1}z^{j_1}a^{i_2}z^{j_2}\ldots a^{i_k})=d_{i_1}\ldots d_{i_k}$$
where $\{d_k\}$ is the sequence of numbers in Definition 9.1.
\end{theorem}

\begin{proof}
First, it follows from definitions that the elements of $N=<aza>$ are indeed as those in the statement, i.e. with $i_1\neq 0$, $i_k\neq 0$.

We define a function $d:N\to\mathbb R$ by the formula in the statement.

{\bf Step 1.} Our first claim is that $d$ is indeed the dimension, for any of the basic corepresentations $r_j=az^ja$. Indeed, for $j\in\{1,\ldots,s-1\}$ we have:
\begin{eqnarray*}
d(az^ja)
&=&d_1^2\\
&=&n\\
&=&\dim (r_j)\\
&=&\dim(az^ja)
\end{eqnarray*}

Also, in the remaining case $j=0$, we have:
\begin{eqnarray*}
d(a^2)
&=&d_2\\
&=&n-1\\
&=&\dim(r_0)\\
&=&\dim(a^2)
\end{eqnarray*}

{\bf Step 2.} Our second claim is that $d$ is a fusion semiring morphism. 

Indeed, we can prove this by using the recurrence relation for the fusion rules on $N$ from Theorem 8.2. It is enough to check that the morphism property holds when applying $d$ to both terms, and we can do this by using a recurrence on $\Sigma i_p$. Indeed, we have: 
\begin{eqnarray*}
d(a^{i_1}z^{j_1}\ldots a^{i_k}\otimes a^{I_1}z^{J_1}\ldots a^{I_K})
&=&d(a^{i_1}z^{j_1}\ldots z^{j_{k-1}}a^{i_k+I_1}z^{J_1}\ldots z^{J_{K-1}}a^{I_K})\\
&&+d(a^{i_1}z^{j_1}\ldots a^{i_k-1}\otimes a^{I_1-1}z^{J_1}\ldots a^{I_K})\\
&=&d_{i_1}\ldots d_{i_{k-1}}d_{i_k+I_1}d_{I_2}\ldots d_{I_K}\\
&&+d_{i_1}\ldots d_{i_{k-1}}d_{i_k-1}\cdot d_{I_1-1}d_{I_2}\ldots d_{I_K}
\end{eqnarray*}

Now by using the Clebsch-Gordan type formula in Proposition 9.2, we can complete the proof of the recurrence step:
\begin{eqnarray*}
d(a^{i_1}z^{j_1}\ldots a^{i_k}\otimes a^{I_1}z^{J_1}\ldots a^{I_K})
&=&d_{i_1}\ldots d_{i_{k-1}}(d_{i_k+I_1}+d_{i_k-1}d_{I_1-1})d_{I_2}\ldots d_{I_K}\\
&=&d_{i_1}\ldots d_{i_{k-1}}(d_{i_k}d_{I_1})d_{I_2}\ldots d_{I_K}\\
&=&(d_{i_1}\ldots d_{i_{k-1}}d_{i_k})(d_{I_1}d_{I_2}\ldots d_{I_K})\\
&=&d(a^{i_1}z^{j_1}\ldots a^{i_k})d(a^{I_1}z^{J_1}\ldots a^{I_K})
\end{eqnarray*}

{\bf Step 3.} We conclude the proof. Since both $d$ and the usual dimension are morphisms of semirings $\mathbb NN\to (\mathbb R^+,+,\cdot)$, which coincide on the basic generators $r_j=az^ja$, these two morphisms are equal, and we are done.
\end{proof}

\section{Concluding remarks}

The fusion rules for the algebra $A_h^s(n)$, computed in this paper, appear to be quite similar to the fusion rules for the algebras $A_o(n),A_u(n)$, previously computed in \cite{ba1}. In what follows we present an attempt of unification.

First, we have the following quite technical definition, from \cite{bbc}.

\begin{definition}
A free quantum algebra is a Hopf algebra $A$ satisfying 
$$A_u(n)\to A\to A_s(n)$$
and having the property that its tensor category is spanned by partitions.
\end{definition}

We should mention that this axiomatization is not fully satisfactory, in the sense that, while being quite restrictive, it still allows too many algebras. The correct axioms, not known so far, should include some conditions which guarantee the compatibility with the Bercovici-Pata bijection \cite{bpa}. See \cite{bb+}.

Now regarding the fusion rules, a possible axiomatization is as follows. Let $R$ be a set, given with two maps as follows:
\begin{enumerate}
\item A map $R\to R$, denoted $r\to\bar{r}$ and called involution.
\item A map $R\times R\to R\cup\{\emptyset\}$, denoted $(r,s)\to r\cdot s$ and called fusion.
\end{enumerate}

The involution and fusion operation extend to the free monoid $<R>$ formed by the words in elements of $R$, in the following way:
$$(r_1\ldots r_k)^-=\bar{r}_k\ldots\bar{r}_1$$
$$(r_1\ldots r_k)\cdot(s_1\ldots s_l)=r_1\ldots r_{k-1}(r_k\cdot s_1)s_2\ldots s_l$$

In the case $r_k\cdot s_1=\emptyset$, or when $k$ or $l$ vanishes, the convention is that the whole word dissapears.

\begin{definition}
A free fusion semiring is a free monoid $<R>$ with fusion rules of the form
$$x\otimes y=\sum_{x=vz,y=\bar{z}w}vw+v\cdot w$$
where $R$ is a set, with involution $r\to\bar{r}$ and fusion $(r,s)\to r\cdot s$. 
\end{definition}

The fusion semirings computed so far are all free, the results in \cite{ba1} and in the present paper being as follows:
\begin{enumerate}
\item For $A_o(n)$ we have $R=\{1\}$, with $\bar{r}=r$ and $r\cdot s=\emptyset$.

\item For $A_u(n)$ we have $R=\mathbb Z_2$, with $\bar{r}=1-r$ and $r\cdot s=\emptyset$.

\item For $A_h^s(n)$ we have $R=\mathbb Z_s$, with $\bar{r}=-r$ and $r\cdot s=r+s$.
\end{enumerate}

We did as well some extra computations, for certain free quantum algebras in preparation, and the conclusion is as follows: (1) the fusion semiring is always free, and (2) the data $(R,-,\cdot)$ doesn't seem to be further axiomatizable.

Summarizing, the question that we would like to raise is as follows.

\begin{conjecture}
If $A$ is free then $R^+(A)$ is free.
\end{conjecture}

It is our hope that further advances on this question will ultimately lead to the technical ingredients needed in order to extend the various analytic results in \cite{ba1}, \cite{bve}, \cite{dvv}, \cite{vvv}, \cite{vve}, \cite{ve1}, \cite{ve2}, \cite{ve3} to the arbitrary free quantum algebras.

We would like to end by recalling an important statement in this direction.

Let $A$ be a Hopf $C^*$-algebra in the sense of Woronowicz. It is known from \cite{wo1} that $A$ has a unique Haar functional, which is not necessarily faithful. By dividing $A$ by the null ideal of the Haar functional we obtain the reduced algebra $A_{red}$. The map $A\to A_{red}$ is an isomorphism when $A$ is amenable in the discrete quantum group sense, and all the standard amenability statements for discrete groups extend to this situation. See B\'edos, Conti and Tuset \cite{bct}.

We recall also that a $C^*$-algebra $A$ has the Dixmier property when any element can be averaged with unitaries, as to get arbitrarily close to the scalars. When $A$ has a trace (and this is the case with $A_{red}$), this is the same as saying that $A$ is simple, and its trace is unique. See Haagerup and Zsido \cite{hzs}.

\begin{conjecture}
If $A$ is free then $A_{red}$ has the Dixmier property.
\end{conjecture}

The first piece of evidence comes from the verification in \cite{ba1} for the Wang algebra $A_u(n)$. The proof there follows Powers' method in \cite{pow}, as modified by de la Harpe and Skandalis in \cite{hsk}, and heavily relies on the fusion rules.

The other piece of evidence comes from the verification in \cite{vve} for the Wang algebra $A_o(n)$. The proof there, which is much more technical, makes use of the notion of Wenzl projection, which ultimately comes from the freeness of $R^+$.

Summarizing, the results in this paper provide a strong evidence for Conjecture 10.3, which in turn would be an important step towards proving Conjecture 10.4.

Let us also mention that a global approach to these problems seems to emerge from the recent work of De Rijdt, Vaes and Vander Vennet \cite{dvv}, \cite{vvv}. So, our first question would be if the methods there can be applied to $A_h^s(n)$.

\end{document}